\documentclass[oneside,english]{amsart}
\usepackage[T1]{fontenc}
\usepackage[latin9]{inputenc}
\usepackage{refstyle}
\usepackage{amsthm}
\usepackage{amssymb}

\makeatletter

\AtBeginDocument{\providecommand\subref[1]{\ref{sub:#1}}}
\RS@ifundefined{subref}
  {\def\RSsubtxt{section~}\newref{sub}{name = \RSsubtxt}}
  {}
\RS@ifundefined{thmref}
  {\def\RSthmtxt{theorem~}\newref{thm}{name = \RSthmtxt}}
  {}
\RS@ifundefined{lemref}
  {\def\RSlemtxt{lemma~}\newref{lem}{name = \RSlemtxt}}
  {}

\numberwithin{equation}{section}
\numberwithin{figure}{section}
\theoremstyle{plain}
\newtheorem{thm}{\protect\theoremname}[section]
  \theoremstyle{remark}
  \newtheorem{rem}[thm]{\protect\remarkname}
  \theoremstyle{remark}
  \newtheorem*{rem*}{\protect\remarkname}

\newref{prop}{name=Proposition~}
\newref{lem}{name=Lemma~}
\newref{thm}{name=Theorem~}
\newref{rem}{name=Remark~}
\newref{cor}{name=Corollary~}

\makeatother

\usepackage{babel}
  \providecommand{\remarkname}{Remark}
\providecommand{\theoremname}{Theorem}

\begin{document}

\title{The real numbers - a survey of constructions}

\author{Ittay Weiss}
\begin{abstract}
We present a comprehensive survey of constructions of the real numbers
(from either the rationals or the integers) in a unified fashion,
thus providing an overview of most (if not all) known constructions
ranging from the earliest attempts to recent results, and allowing
for a simple comparison-at-a-glance between different constructions. 
\end{abstract}

\maketitle
\begin{tiny}\tableofcontents{}\end{tiny}

\section{Introduction}

The novice, through the standard elementary mathematics indoctrination,
may fail to appreciate that, compared to the natural, integer, and
rational numbers, there is nothing simple about defining the real
numbers. The gap, both conceptual and technical, that one must cross
when passing from the former to the latter is substantial and perhaps
best witnessed by history. The existence of line segments whose length
can not be measured by any rational number is well-known to have been
discovered many centuries ago (though the precise details are unknown).
The simple problem of rigorously introducing mathematical entities
that do suffice to measure the length of any line segment proved very
challenging. Even relatively modern attempts due to such prominent
figures as Bolzano, Hamilton, and Weierstrass were only partially
rigorous and it was only with the work of Cantor and Dedekind in the
early part of the 1870's that the reals finally came into existence.
The interested reader may consult \cite{ebbinghaus1990numbers} for
more on the historical developments and further details. 

Two of the most famous constructions of the reals are Cantor's construction
by means of Cauchy sequences of rational numbers and Dedekind's construction
by means of cuts of rational numbers, named after him. Detailed accounts
of these constructions and their ubiquity in textbooks, together with
the well-known categoricity of the axioms of a complete ordered field,
would have put an end to the quest for other constructions, and yet
two phenomena persist. Firstly, it appears that human curiosity concerning
the real numbers is not quite quenched with just these two constructions.
Even though any two models of the axioms of a complete ordered field
are isomorphic, so it really does not matter which model one works
with, we still seem to be fascinated with finding more and more different
models to the same abstract concept. Secondly, and more practically,
from the constructive point of view not all models of the real numbers
are isomorphic. Fueled by applications in automated theorem proving
and verification, where one must represent the real numbers in a computer,
nuances of the differences between various constructions of the reals
become very pronounce. We refer the reader to \cite{bridges1994constructive}
and \cite{harrison1994constructing,harrison1998theorem} for more
details on the constructive reals and on theorem proving with the
real numbers, respectively.

\section{A survey of constructions\label{sec:A-survey-of}}

In order to present a uniform survey of constructions of the real
numbers we choose to adopt the following somewhat debatable point
of view according to which every construction of the real numbers
ultimately relies on an \emph{observation} about the reals (treated
axiomatically) leading to a bijective correspondence between the set
of real numbers and a set defined in terms of simpler entities (often
the rational or the integer numbers) upon which agreement of existence
is present. That set is then taken to be the definition of the reals,
with the \emph{order }structure and the \emph{arithmetical }operations
defined, examined, and eventually shown to form a complete ordered
field. 

We present what we hope is an exhaustive list of constructions of
the reals one can find in the literature, all following the presentation
style exposited above. Certainly, this restrictive decision sometimes
necessitates a suboptimal presentation of a particular construction,
however the uniform style makes comparison between the definitions
easier. 

As a convention, let $\mathbb{N}^{+}=\mathbb{N}\cup\{\omega\}$ be
the set of all natural numbers augmented with the symbol $\omega$
which algebraically behaves like $\infty$. In particular, $x\le\omega$
for all $x\in\mathbb{N}^{+}$ and we define $x+\omega=\omega=\omega+x$
and $x\omega=\omega=\omega x$ for all $x\in\mathbb{N}^{+}$, and
$\frac{x}{\omega}=0$ for all $x\in\mathbb{Z}$. The sole use of this
convention is in treating finite sequences of integers as infinite
ones ending with a constant stream of $\omega$'s. 

Finally, we mention at this point, rather than at each construction
surveyed below, that typically it makes little difference whether
one constructs the positive (or nonnegative) reals $\mathbb{R}_{+}$
and then extend to all the reals by formally adding inverses (and
a $0$ if needed), or constructing all of $\mathbb{R}$ in one go.
However, the former approach may sometimes be technically simpler
than the latter. Consequently, below, a survey of a construction will
be considered complete even if it only produces $\mathbb{R}_{+}$. 
\begin{rem}
In those constructions below that refer to convergence of a sequence
$(q_{n})$ of rationals to a rational number $q$, the precise meaning
of such a statement is that for every rational number $\varepsilon>0$
there exists $n_{0}\in\mathbb{N}$ with $|q_{n}-q|<\varepsilon$ for
all $n>n_{0}$. A sequence that converges to $0$ is also known as
a \emph{null sequence}. 
\end{rem}

\subsection{Stevin's ``construction'' (\emph{De Thiende} and \emph{L'Arithmetique},
1585)\label{sub:Stevin}}
\begin{rem*}
Stevin is credited with laying down the foundations of the the decimal
notation. Stevin did not produce a rigorous construction of the reals
though he did present the then controversial point of view that there
is nothing significantly different in nature between the rational
numbers and the irrational ones. Constructing the reals as decimal
expansions (or in any other base) is a popular approach by novices
but is fraught with technical difficulties. The allure of this approach
most likely lies in the emphasis decimal expansions receive in the
current mathematical curriculum, where decimal expansions triumph
over anything else. The details presented below are nowhere near what
Stevin presented. Instead, we follow Gowers (\cite{Gowers}) but leave
the details at a minimum. We also mention \cite{katz2012stevin} for
a more holistic view of Stevin's numbers. 
\end{rem*}

\subsubsection*{Observation}

Every real number $a$ can be written as 
\[
a=\sum_{k=-\infty}^{\infty}\frac{a_{k}}{b^{k}}
\]
where the \emph{base }$b\ge2$ is an integer and the $a_{k}\in\{0,1,2,\dots,b-1\}$
are digits, and there exists $k_{0}\in\mathbb{Z}$ such that $a_{k}=0$
for all $k<k_{0}$. Moreover, the presentation is unique if we further
demand that there does not exist $k_{0}\in\mathbb{Z}$ with $a_{k}=b-1$
for all $k>k_{0}$.

\subsubsection*{The reals}

One may now take the formal $b$-base expansions as above to be the
real numbers.

\subsubsection*{Order}

Defining the order between $b$-base expansions presents no difficulties;
$a<a'$ precisely when $a_{k_{0}}<a_{k_{0}}'$ for the largest index
$k_{0}$ with $a_{k_{0}}\ne a'_{k_{0}}$.

\subsubsection*{Arithmetic}

The algorithms for symbolically performing addition and multiplication
of real numbers are cumbersome. Gowers suggests that the simplest
approach to turn Stevin's $b$-base expansions into a construction
of the reals is by employing limiting arguments to define addition
and multiplication.

\subsection{Weierstrass's construction (unpublished by Weierstrass, we paraphrase
following \cite{tweddle2011weierstrass}, ca. 1860)\label{sub:Weierstrass}}

\subsubsection*{Observation}

Every positive real number $a$ can be written as 
\[
a=\sum_{s\in S}s
\]
where $S$ is a multiset (i.e., a set where elements may be repeated
more than once) whose elements consist of positive integers and positive
rationals of the form $\frac{1}{n}$. Such a presentation is, of course,
not unique.

\subsubsection*{The reals}

Consider the set $\mathbb{S}$ of all non-empty multisets $S$ of
positive integers and positive rationals of the form $\frac{1}{n}$,
which are \emph{bounded }in the sense that there exists $M>0$ with
\[
\sum_{s\in S_{0}}s<M
\]
 for all finite submultisets $S_{0}\subseteq S$ (being finite means
that the total number of elements in $S_{0}$, counting multiplicities,
is finite). Declare for two multisets $S,T\in\mathbb{S}$ that $S\le T$
if for every finite submultiset $S_{0}\subseteq S$ there exists a
finite submultiset $T_{0}\subseteq T$ with
\[
\sum_{s\in S_{0}}s\le\sum_{t\in T_{0}}t.
\]
Declare $S\sim T$ if both $S\le T$ and $T\le S$ hold. The set $\mathbb{R}_{+}$
of positive real numbers is then defined to be $\mathbb{S}/{\sim}$.

\subsubsection*{Order}

For two real numbers $a=[S]$ and $b=[T]$, the relation $a\le b$
holds when $S\le T$.

\subsubsection*{Arithmetic}

Addition and multiplication of positive integers and of positive rationals
of the form $\frac{1}{n}$ extends to $\mathbb{S}$ by $S+T=\{s+t\mid s\in S,t\in T\}$
and $ST=\{st\mid s\in S,t\in T\}$ subject to the convention that
multiplicities are taken into consideration and that any sum or product
which is not of the form of an integer or $\frac{1}{n}$ is replaced
by a (necessarily) finite number of elements of this form (there is
no canonical choice, and any would do). 

The addition and multiplication of the real numbers $a=[S]$ and $b=[T]$
is given by
\[
a+b=[S+T]
\]
and 
\[
ab=[ST].
\]

\begin{rem}
As noted, the construction lacks in rigor. 
\end{rem}

\subsection{Dedekind's construction (\cite{dedekind1930stetigkeit}, 1872)\label{sub:Dedekind}}

\subsubsection*{Observation}

Any real number $a$ determines a partition of $\mathbb{Q}$ into
a pair $(A,B)$ where $A=\{q\in\mathbb{Q}\mid q<a\}$ and $B=\{q\in\mathbb{Q}\mid a\le q\}$.
Obviously, $A$ is non-empty and downward closed, $B$ is non-empty
and upward closed, and $A$ has no greatest element. Any partition
of $\mathbb{Q}$ satisfying these properties is called a \emph{Dedekind
cut} and this construction is a bijection between the real numbers
and the set of all cuts.

\subsubsection*{The reals}

The set $\mathbb{R}$ of real numbers is defined to be the set of
all Dedekind cuts.

\subsubsection*{Order}

Given real numbers $a_{1}=(A_{1},B_{1})$ and $a_{2}=(A_{2},B_{2})$,
the relation $a_{1}\le a_{2}$ holds precisely when $A_{1}\subseteq A_{2}$.

\subsubsection*{Arithmetic}

Obviously, in any Dedekind cut $(A,B)$, any one of $A$ or $B$ determines
the other and if $A\subsetneq\mathbb{Q}$ satisfies the properties
of left 'half' of a Dedekind cut, then $(A,\mathbb{Q}\setminus A)$
is a Dedekind cut. It thus suffices to concentrate on $A$. Addition
of real numbers given by Dedekind cuts represented by sets $A_{1}$
and $A_{2}$ is defined by 
\[
A_{1}+A_{2}=\{a_{1}+a_{2}\mid a_{1}\in A_{1},a_{2}\in A_{2}\}.
\]
If $A_{1}$ and $A_{2}$ represent non-negative reals, then their
product is given by 
\[
A_{1}A_{2}=\{q\in\mathbb{Q}\mid q\le0\}\cup\{a_{1}a_{2}\mid a_{1}\in A_{1},a_{2}\in A_{2},a_{1}\ge0,a_{2}\ge0\}.
\]
Multiplication is then extended by sign cases as usual.

\subsection{Cantor's construction (\cite{cantor1966grundlagen}, 1873)\label{sub:Cantor}}

\subsubsection*{Observation}

Every real number $a$ is the limit of a sequence $(q_{n})$ of rationals.
Moreover, any two convergent sequences $(q_{n})$ and $(q_{n}')$
converge to the same value $a$ if, and only if, $|q_{n}-q'_{n}|\xrightarrow[n\to\infty]{}0$.

\subsubsection*{The reals}

Declare a sequence $(q_{n})$ of rationals to be a \emph{Cauchy sequence
}if for all $\varepsilon>0$ there exists $k_{0}\in\mathbb{N}$ with
$|q_{n}-q_{m}|<\varepsilon$, provided that $n,m>k_{0}$. Let $C$
be the set of all Cauchy sequences of rational numbers and let $\sim$
be the equivalence relation on $C$ given by $(q_{n})\sim(q_{n}')$
precisely when $|q_{n}-q'_{n}|\to0$. The set of real numbers is then
$\mathbb{R}=C/{\sim}$, the set of all equivalence classes of Cauchy
sequences modulo $\sim$.

\subsubsection*{Order}

Declare that two real numbers $a=[(q_{n})]$ and $b=[(q_{n}')]$ satisfy
$a<b$ when $a\ne b$ and when there exists $k_{0}\in\mathbb{N}$
with $q_{n}<q_{n}'$ for all $n>k_{0}$.

\subsubsection*{Arithmetic}

Addition and multiplication are given by $a+b=[(q_{n}+q_{n}')]$ and
$ab=[(q_{n}q_{n}')]$, respectively.

\subsection{Bachmann's construction (\cite{Bachmann1892}, 1892)\label{sub:Bachmann}}

The details below are essentially identical to those given by Bachmann,
but the style is slightly adapted.

\subsubsection*{Observation}

A sequence $\{I_{n}\}_{n\ge1}$ of intervals $I_{n}=[a_{n},c_{n}]$
in the real line is said to be a \emph{nested family of intervals
}or more simply a \emph{nest}, if $I_{k+1}\subseteq I_{k}$ for all
$k\ge1$ and $c_{n}-a_{n}\xrightarrow[n\to\infty]{}0$. For each such
nest there is then a unique real number $b$ satisfying $b\in I_{k}$
for all $k\ge1$. Moreover, two nests determine in this way the same
real number if, and only if, the nests admit a common refinement.
In more detail, a nest $\{I_{n}\}$ is \emph{finer} than a nest $\{J_{n}\}$
when $I_{n}\subseteq J_{n}$, for all $n\ge1$. Two nests have a common
refinement if there is a nest finer than each of them. Due to the
density of the rational numbers in the real numbers the intervals
above can be replaced by rational intervals consisting of rational
numbers only, while retaining the correspondence with the reals.

\subsubsection*{The reals}

Consider now \emph{rational intervals} of the form $I=[a,c]=\{x\in\mathbb{Q}\mid a\le x\le c\}$,
where $a,c\in\mathbb{Q}$. A \emph{rational nest} is a family $\{I_{n}\}_{n\ge1}$
of rational intervals $I_{n}=[a_{n},c_{n}]$ satisfying $I_{k+1}\subseteq I_{k}$
for all $k\ge1$ and $c_{n}-a_{n}\xrightarrow[n\to\infty]{}0$. A
rational nest $\{I_{n}\}$ is \emph{finer }than a rational nest $\{J_{n}\}$
if $I_{n}\subseteq J_{n}$ for all $n\ge1$. Consider now the set
$N$ of all rational nests, and define on it the relation $\sim$
whereby $\{I_{n}\}\sim\{I'_{n}\}$ precisely when there exists a common
refinement of $\{I_{n}\}$ and $\{I'_{n}\}$. It follows easily that
$\sim$ is an equivalence relation on $N$ and the set $\mathbb{R}$
of real numbers is defined to be $N/{\sim}$, the set of equivalence
classes of rational nests.

\subsubsection*{Order}

Two real numbers $x=[\{I_{n}\}]$ and $y=[\{J_{n}\}]$ satisfy $x<y$
precisely when there exists $n_{0}\in\mathbb{N}$ with $I_{n_{0}}<J_{n_{0}}$
in the sense that $\alpha<\beta$ for all $\alpha\in I_{n_{0}}$ and
all $\beta\in J_{n_{0}}$.

\subsubsection*{Arithmetic}

Extending the arithmetic operations of addition and multiplication
of rational numbers to subsets $S,T$ of rational numbers by means
of $S+T=\{s+t\mid s\in S,t\in T\}$ and $ST=\{st\mid s\in S,t\in T\}$,
it is easily seen that for all rational intervals $I$ and $J$, both
$I+J$ and $IJ$ are again rational intervals. Addition and multiplication
of the real numbers $x$ and $y$ is given by $x+y=[\{I_{n}+J_{j}\}]$
and $xy=[\{I_{n}J_{n}\}]$.

\subsection{Bourbaki's approach to the reals (\cite{bourbaki1998general}, ca.
1960)}

Bourbaki develops the general machinery of uniform spaces and their
completion, observes that the rationals admit a uniform structure,
and takes $\mathbb{R}$ to be any completion of $\mathbb{Q}$. The
structure of $\mathbb{R}$ as a complete ordered field is then deduced
using the machinery of uniform spaces. Strictly speaking then, Bourbaki
does not construct the reals, and in fact stresses the point that
no particular construction is required; the universal properties,
provided by any completion, suffice. However, Bourbaki also (famously)
discusses a particular completion process of any uniform space (initially
by means of equivalence classes of Cauchy filters, with the canonical
choice of minimal Cauchy filters later on). Bourbaki's constructions
can be combined and expanded into a particular construction of the
reals, which we thus refer to as the Bourbaki reals, and we cast them
into the mold of the survey. The details of this construction will
be given in a separate article.

\subsubsection*{Observation}

Any real number $x$ gives rise to two filters, namely the principal
filter $\langle x\rangle=\{S\subseteq\mathbb{R}\mid x\in S\}$ and
the minimal Cauchy filter $\iota(x)$ generated by all intervals containing
$x$. Each construction leads to a bijective correspondence between
$\mathbb{R}$ and a certain set of filters on $\mathbb{R}$, but the
latter can be used to obtain a bijection between $\mathbb{R}$ and
the set of all minimal Cauchy filters on $\mathbb{Q}$ by means of
simply intersecting each set in $\iota(x)$ with $\mathbb{Q}$.

\subsubsection*{The reals}

The set $\mathbb{R}$ of real numbers is defined to be the set of
all minimal Cauchy filters on $\mathbb{Q}$.

\subsubsection*{Order}

The order relation on $\mathbb{Q}$ extends universally to $\mathcal{P}(\mathbb{Q})$
by declaring $A<_{\forall}B$ when for all $a\in A$ and $b\in B$
one has $a<b$. Similarly, the relation $<_{\forall}$ on $\mathcal{P}(\mathbb{Q})$
extends existentially to $\mathcal{P}(\mathcal{P}(\mathbb{Q}))$ by
declaring $\mathcal{A}<_{\exists\forall}\mathcal{B}$ when there exist
$A\in\mathcal{A}$ and $B\in\mathcal{B}$ with $A<_{\forall}B$. None
of these extensions is an order relation. However, since any filter
on $\mathbb{Q}$ is an element in $\mathcal{P}(\mathcal{P}(\mathbb{Q}))$,
the relation $<_{\exists\forall}$ restrict to a relation on $\mathbb{R}$,
and this relation is an order relation.

\subsubsection*{Arithmetic}

Addition and multiplication in $\mathbb{Q}$ extend element-wise to
subsets $A,B\subseteq\mathbb{Q}$. Given real numbers $x$ and $y$,
i.e., minimal Cauchy filters on $\mathbb{Q}$, their sum is $x+y=\langle\{A+B\mid A\in a,B\in b\}\rangle$
and their product is $xy=\langle\{AB\mid A\in a,B\in b\}\rangle$.
It should be noted that the fact that $x+y$ and $xy$ are real numbers,
i.e., that the defining generated filters are \emph{minimal}, is not
a triviality but rather a fact that encapsulates most of the technical
details in the construction, rendering the rest of the proof quite
straightforward.
\begin{rem}
This construction can be seen as Bachmann's construction (\subref{Bachmann})
with a canonical choice of representatives. A direct comparison from
Bachmann's reals to the Bourbaki reals is given by sending a nest
of intervals to the roundification of the filter generated by the
intervals. 
\end{rem}

\subsection{Maier-Maier's construction by a variation on Dedekind cuts (\cite{maier1973},
1973)\label{sub:Maier-Maier}}

\subsubsection*{Observation}

Every real number $a$ occurs as the greatest lower bound of the set
$\{q\in\mathbb{Q}\mid a<q\}$. Of course, the same real number is
the greatest lower bound of many other subsets of $\mathbb{Q}$. However,
two bounded below sets $T_{1},T_{2}\subseteq\mathbb{Q}$ have the
same greatest lower bound provided that the set of lower bounds of
$T_{1}$ coincides with the set of lower bounds of $T_{2}$.

\subsubsection*{The reals}

Let $B$ be the set of all subsets of $\mathbb{Q}$ which are bounded
below, and denote, for $T\in B$, by $b(T)$ the set of all lower
bounds of $T$. Given $T_{1},T_{2}\in B$, declare that $T_{1}\sim T_{2}$
precisely when $b(T_{1})=b(T_{2})$. It is easily seen that $\sim$
is an equivalence relation, and the real numbers are defined to be
$B/{\sim}$, the set of equivalence classes.

\subsubsection*{Order}

Given real numbers $x=[S]$ and $y=[T]$, the relation $x<y$ holds
precisely when $b(S)\subset b(T)$.

\subsubsection*{Arithmetic}

For real numbers $x$ and $y$, their sum is given by $x+y=[\{s+t\mid s\in S,t\in T\}]$.
The product of $x$ and $y$, provided that all the elements in $S$
and in $T$ are positive, is given by $xy=[\{st\mid s\in S,t\in T\}]$.
Multiplication is extended to all real numbers by sign considerations. 
\begin{rem}
This construction is essentially Dedekind's construction without canonical
choices of representatives. In more detail, given a real number $x=[S]$
the set $b(S)$, if it does not have a maximum determines a Dedekind
cut, denoted by $b_{0}(x)$. If $b(S)$ does have a maximum, then
$b(S)\setminus\{\max b(S)\}$ determines a Dedekind cut, again denoted
by $b_{0}(x)$. The function $x\mapsto b_{0}(x)$ is then an isomorphism
giving a direct comparison between Dedekind's construction and the
present construction. 
\end{rem}

\subsection{Shiu's construction by infinite series (\cite{shiu1974}, 1974)\label{sub:Shiu}}

\subsubsection*{Observation}

Since the harmonic series $\sum_{n}\frac{1}{n}$ diverges, it follows
that every positive real number $x$ can be written as 
\[
x=\sum_{n\in A}\frac{1}{n}
\]
for some (non-unique) infinite set $A\subseteq\mathbb{N}$.

\subsubsection*{The reals}

Let $\alpha$ be the set of all infinite subsets of natural numbers,
writing $A=(a_{k})$ with $a_{k}<a_{k+1}$ for a typical element in
$\alpha$. For such an $A\in\alpha$ let 
\[
A_{n}=\sum_{k=1}^{n}\frac{1}{a_{k}}.
\]
Let $\beta$ be the subset of $\alpha$ consisting of those $A\in\alpha$
for which the sequence $(A_{n})$ is bounded. Introduce an equivalence
relation on $\beta$ by declaring $A\sim B$ precisely when $(A_{n}-B_{n})$
is a null sequence. The positive real numbers are then defined to
be $\mathbb{R}_{+}=\beta/{\sim}$, the set of equivalence classes.

\subsubsection*{Order}

Real numbers $x$ and $y$ satisfy $x\le y$ precisely when $x=[A]$
and $y=[B]$ for some representatives satisfying $A\subseteq B$.

\subsubsection*{Arithmetic}

Let $x=[A]$ and $y=[B]$ be positive real numbers. Consider the set
$AB=\{ab\mid a\in A,b\in B\}$. Call the representing sets $A$ and
$B$ \emph{excellent }if $A\cap B=\emptyset$ and every $c\in AB$
can be written uniquely as $c=ab$ with $a\in A$ and $b\in B$. Heuristically,
\[
x=\sum_{a\in A}\frac{1}{a}
\]
and 
\[
y=\sum_{b\in B}\frac{1}{b}
\]
so, since the representatives are excellent, it follows that both
$A\cup B$ and $AB$ represent real numbers, which intuitively are
\[
\sum_{c\in A\cup B}\frac{1}{c}=\sum_{a\in A}\frac{1}{a}+\sum_{b\in B}\frac{1}{b}=x+y
\]
and 
\[
\sum_{c\in AB}\frac{1}{c}=(\sum_{a\in A}\frac{1}{a})(\sum_{b\in B}\frac{1}{b})=xy.
\]
This informal argument turns into a definition of addition and multiplication
on representatives by the fact that excellent representatives can
always be found. 
\begin{rem}
With suitable adaptation the harmonic series can be replaced by other
divergent series of positive rationals converging to $0$. 
\end{rem}

\begin{rem}
This construction is very similar to Weierstrass's (\subref{Weierstrass}).
Here repetitions are not allowed, resulting in a simpler definition
of the set of real numbers at the cost of a slightly less immediate
notion of addition and multiplication. Weierstrass allows repetitions
and thus arithmetic is immediate, but identifying the set of real
numbers requires a delicate notion of equivalence. 
\end{rem}

\subsection{Faltin-Metropolis-Ross-Rota's wreath construction (\cite{faltin1975},
1975)\label{sub:Faltin}}

\subsubsection*{Observation}

The difficulty in defining the arithmetic operations when defining
the reals as sequences of base $b$ expansions lies in the need to
keep track of carries. This necessity stems from the (almost) uniqueness
of the digits of any given real number, resulting from the use of
the base to limit the range of the digits. Instead one may not place
a limit on the digits, i.e., every real number $a$ can be written
in infinitely many ways as 
\[
a=\sum_{n\in\mathbb{Z}}\frac{a_{n}}{2^{n}}
\]
where the $a_{n}$ are integers, all of which are $0$ for sufficiently
small $n$ (the base is taken to be $b=2$ only to conform with the
construction in the mentioned article). The definition of addition
and multiplication of such expansions is formally identical to the
way one would add and multiply formal Laurent series, at the price
of an algorithmically more intricate recovery of the order structure
by manipulating carries. Interestingly, this exchange in algorithmic
complexity between arithmetic and order results in a much simpler
construction of the real numbers than Stevin's construction.

\subsubsection*{The reals}

Let $\Sigma(\mathbb{Z})$ be the set of all formal expressions of
the form
\[
\sum_{n\in\mathbb{Z}}a_{n}x^{n}
\]
where $x$ is an indeterminate, $a_{n}\in\mathbb{Z}$ for all $n\in\mathbb{Z}$,
and $a_{k}=0$ for all $k<k_{0}$, for some $k_{0}\in\mathbb{Z}$.
With the formal operations of addition and multiplication of Laurent
series the set $\Sigma(\mathbb{Z})$ becomes a ring, whose elements
are also called \emph{strings}. More explicitly, 
\[
\sum_{n\in\mathbb{Z}}a_{n}x^{n}+\sum_{n\in\mathbb{Z}}b_{n}x^{n}=\sum_{n\in\mathbb{Z}}(a_{n}+b_{n})x^{n}
\]
and
\[
(\sum_{n\in\mathbb{Z}}a_{n}x^{n})(\sum_{n\in\mathbb{Z}}b_{n}x^{n})=\sum_{n\in\mathbb{Z}}c_{n}x^{n}
\]
where 
\[
c_{k}=\sum_{n\in\mathbb{Z}}a_{n}b_{k-n},
\]
for all $k\in\mathbb{Z}$. The element $K\in\Sigma(\mathbb{Z})$ with
$k_{0}=1$, $k_{1}=-2$, and $k_{i}=0$ for all $i\in\mathbb{Z}\setminus\{0,1\}$
is called the \emph{carry constant}. It is easily seen that two elements
$A,B\in\Sigma(\mathbb{Z})$ are related by $A=B+KC$, where $C\in\Sigma(\mathbb{Z})$
has only finitely many non-zero coefficients, precisely when $A$
can be obtained from $B$ by formally performing carrying operations
as indicated by $C$ (in base $2$). An element $A\in\Sigma(\mathbb{Z})$
is said to be \emph{bounded }if there exists an integer $z\ge1$ such
that 
\[
\sum_{i\le n}|a_{i}|2^{n-i}\le z2^{n}.
\]
for all non-negative $n$. The set of all bounded elements in $\Sigma(\mathbb{Z})$
is denoted by $\Sigma_{2}(\mathbb{Z})$. An element $C\in\Sigma(\mathbb{Z})$
is called a \emph{carry string }if $KC$ is bounded, and when for
every positive integer $z$ there exists $k\ge0$ with 
\[
z|c_{j}|\le2^{j}
\]
for all $j>k$. Finally, two bounded elements $A,B\in\Sigma_{2}(\mathbb{Z})$
are declared to be equivalent if there exists a carry string $C$
with $A=B+KC$. The set $\mathbb{R}$ of real numbers is then $\Sigma_{2}(\mathbb{Z})/{\sim}$,
the set of equivalence classes of formal carry-free binary expansions
modulo the performance of carrying.

\subsubsection*{Arithmetic}

The ring structure on $\Sigma(\mathbb{Z})$ restricts to one on $\Sigma_{2}(\mathbb{Z})$
and is compatible with the equivalence relation $A=B+KC$, and thus
gives rise to addition and multiplication in $\mathbb{R}$, namely
the usual addition and multiplication of formal Laurent series performed
on representatives.

\subsubsection*{Order}

For the definition of a \emph{clear string }refer to \cite[Section 6]{faltin1975}.
For every real number $[A]$ there exists a unique clear string $B$
with $[A]=[B]$. Declare $[A]$ to be \emph{positive} if, when cleared,
the leading digit (i.e., leading non-zero coefficient) is $1$. Then
the set of positive reals defines an ordering in the usual manner.
Equivalently, $[A]\le[B]$ is the lexicographic order on the cleared
strings representing $[A]$ and $[B]$.

\subsection{De Bruijn's construction by additive expansions (\cite{debruijn1976},
1976)\label{sub:De-Bruijn}}

\subsubsection*{Observation}

As noted in \subref{Stevin}, the set of real numbers can be identified
with formal decimal expansions (or other bases), i.e., as certain
strings of digits indexed by the integers. The difficulty of performing
the arithmetical operations (and even just addition) directly on the
strings of digits stems in some sense from the expansions arising
in complete disregard of the arithmetical operations; the expansions
are analytic, not algebraic. If, instead, one considered the set of
formal expansions with the aim of focusing on easily defining addition,
then one is led to interpret the expansions differently. This is the
approach taken in this construction.

\subsubsection*{The reals}

Fix an integer $b>1$ and let $\Sigma$ be the set of all functions
$f:\mathbb{Z}\to\{0,1,2,\dots b-1\}$ which satisfy the condition
that for all $i\in\mathbb{Z}$ there exists $k\in\mathbb{Z}$ with
$k>i$ and $f(k)<b-1$. Given any two functions $f,g\colon\mathbb{Z}\to\{0,1,2,\ldots,b-1\}$
define two other such functions, denoted by ${\rm difcar}(f,g)$ (standing
for the difference carry of $f$ and $g$) and $f-g$, as follows.
For $k\in\mathbb{Z}$ define ${\rm difcar}(f,g)(k)=1$ if there exists
$x\in\mathbb{Z}$ with $x>k$, $f(x)<g(x)$, and such that $f(y)\le g(y)$
for all $k<y<x$. In all other cases ${\rm difcar}(f,g)(k)=0$. The
value of $f-g$ at $k\in\mathbb{Z}$ is given by 
\[
(f-g)(k)=f(x)-g(x)-{\rm difcar}(f,g),
\]
computed $\mod b$. Following this procedure leads to defining $f\in\Sigma$
to be positive if $f\ne0$ and if some $k\in\mathbb{Z}$ exists with
$f(y)=0$ for all $y<k$. Similarly, declare $f\in\Sigma$ to be negative
if there exists $k\in\mathbb{Z}$ with $f(y)=b-1$ for all $y<k$.
Then the set of real numbers is defined to be the set of all $f\in\Sigma$
such that either $f=0$, $f$ is positive, or $f$ is negative.

\subsubsection*{Order}

For real numbers $f$ and $g$, the relation $f<g$ holds precisely
when $g-f$ is positive. The greatest lower bound property is then
verified, allowing for limit-like arguments used only when defining
the product of real numbers.

\subsubsection*{Arithmetic}

Addition of real numbers is given by $f+g=f-(0-g)$. Multiplication
is defined as a supremum over suitably constructed approximations.

\subsection{Rieger's construction by continued fractions (\cite{rieger1982new},
1982)\label{sub:Rieger}}

\subsubsection*{Observation}

Every irrational real number $a$ can be written uniquely as a continued
fraction
\[
a=a_{0}+\frac{1}{a_{1}+\frac{1}{a_{2}+\frac{1}{a_{3}+\cdots}}}=[a_{0};a_{1},a_{2},\dots,a_{k},\ldots]
\]
where $a_{0}\in\mathbb{Z}$ and $a_{k}\in\mathbb{N}$ with $a_{k}\ge1$
for all $k\ge1$. When $a$ is rational the continued fraction terminates
at some $k_{0}\ge0$, and if one further demands that if $k_{0}>0$,
then $a_{k_{0}}>1$, then the presentation of rational numbers is
also unique.

\subsubsection*{The reals}

Let $\mathbb{R}$ be the set of all sequences $[a_{0};a_{1},a_{2},\ldots,a_{k},\dots]$
where $a_{0}\in\mathbb{Z}$ and $a_{k}\in\mathbb{N}^{+}$ with $a_{k}\ge1$
for all $k\ge1$, subject to the demand that if $a_{k}=\omega$, then
$a_{t}=\omega$ for all $t>k$ and if $k_{0}$ is the last index where
$a_{k_{0}}\ne\omega$ and $k_{0}>0$, then $a_{k_{0}}>1$.

\subsubsection*{Order}

Given real numbers $a=[a_{0};a_{1},a_{2},\ldots,a_{k},\ldots]$ and
$b=[b_{0};b_{1},b_{2},\ldots,b_{k},\ldots]$ the relation $a<b$ holds
precisely when $a\ne b$ and when for the smallest index $k_{0}$
with $a_{k_{0}}\ne b_{k_{0}}$ one has
\begin{itemize}
\item $a_{k_{0}}<b_{k_{0}}$, if $k$ is even;
\item $a_{k_{0}}>b_{k_{0}}$, if $k$ is odd.
\end{itemize}
The least upper bound property of $\mathbb{R}$ is then established
and the proof of the Euclidean algorithm produces an order embedding
$\mathbb{Q}\to\mathbb{R}$, which thus serves to identify the rationals
in $\mathbb{R}$ as precisely those real numbers in which $\omega$
appears. It then follows that every real number $a=[a_{0};a_{1},a_{2},\ldots,a_{k},\ldots]$
can be approximated by suitably constructed rationals to obtain
\[
a^{(0)}<a^{(2)}<a^{(4)}<\cdots<a<\cdots<a^{(5)}<a^{(3)}<a^{(1)}.
\]

\subsubsection*{Arithmetic}

The sum of $a$ and $b$ is defined to be 
\[
a+b=\sup\{a^{(2n)}+b^{(2n)}\mid n\ge0\}.
\]
Multiplication of positive real numbers is given by 
\[
ab=\sup\{a^{(2n)}b^{(2n)}\mid n\ge0\}
\]
and extended to all of $\mathbb{R}$ by the usual sign conventions.
The proofs of the algebraic properties utilize the rational approximations
using limit-like arguments.

\subsection{Schanuel (et al)'s construction using approximate endomorphisms of
$\mathbb{Z}$ (\cite{street1985efficient,street2003update,douglas2004efficient,arthan2004,grundhofer2005describing},
1985)}

\subsubsection*{Observation}

Given a real number $a$, the function $f_{a}\colon\mathbb{R}\to\mathbb{R}$
given by $f_{a}(x)=ax$ is a linear function whose slope is $a$,
and the assignment $a\mapsto f_{a}$ thus sets up a bijection between
the real numbers and linear operators $\mathbb{R}\to\mathbb{R}$.
Under this bijection, addition in $\mathbb{R}$ corresponds to the
point-wise addition of functions, while multiplication in $\mathbb{R}$
corresponds to composition of functions. 

Of course, this point of view of the real numbers as linear operators
(thought of as slopes) requires the existence of the real numbers
for the operators to operate on. Thus, in order to obtain a construction
of the reals one seeks to modify $f_{a}$ to a linear operator on
$\mathbb{Z}$ instead of on $\mathbb{R}$. Restricting the domain
of $f_{a}$ to $\mathbb{Z}$ does not produce a function $f_{a}\colon\mathbb{Z}\to\mathbb{Z}$
(unless $a$ is an integer), and it is tempting to simply adjust $f_{a}(x)$
to an integer near $ax$ so as to obtain a function $g_{a}\colon\mathbb{Z}\to\mathbb{Z}$.
Of course, this new function need not be linear any more, but it is
approximately so in the sense that though the choice of $g_{a}$ may
be somewhat arbitrary, as long as the adjustment to an integer was
not too out of hand, the set 
\[
\{g_{a}(x+y)-g_{a}(x)-g_{a}(y)\mid x,y\in\mathbb{Z}\},
\]
which measures how non-linear $g_{a}$ is, is finite. Furthermore,
while it is obvious that different $g_{a},g_{a}'$ may arise from
the same $f_{a}$, for sensible processes leading to a $g_{a}$ and
$g_{a}'$, the set
\[
\{g_{a}(x)-g_{a}'(x)\mid x\in\mathbb{Z}\},
\]
measuring how different $g_{a}$ and $g_{a}'$ are, is finite. Further
motivation is gathered from the observation that defining $g_{a}(x)=[ax]$,
the integer part of $ax$, is a function with the property that $\frac{g_{a}(x)}{x}\xrightarrow[x\to\infty]{}a$,
so there is at least one obvious way of adjusting the linear function
$f_{a}$ to an approximately linear function from which $a$ can be
extracted.

\subsubsection*{The reals}

Let $\mathbb{Z}$ be the integers considered as a group under addition.
Call a function $f\colon\mathbb{Z}\to\mathbb{Z}$ a \emph{quasihomomorphism}
if the set 
\[
\{f(x+y)-f(x)-f(y)\mid x,y\in\mathbb{Z}\}
\]
is finite. Introduce an equivalence relation on the set $H$ of all
quasihomomorphisms whereby $f\sim g$ precisely when the set 
\[
\{f(x)-g(x)\mid x\in\mathbb{Z}\}
\]
is finite. The real numbers are then defined to be $H/{\sim}$.

\subsubsection*{Arithmetic}

Given real numbers $a=[f]$ and $b=[g]$, their sum is represented
by $f+g\colon\mathbb{Z}\to\mathbb{Z}$, where $(f+g)(x)=f(x)+g(x)$.
The product $ab$ is represented by $f\circ g\colon\mathbb{Z}\to\mathbb{Z}$,
the composition of $f$ and $g$.

\subsubsection*{Order}

It can be shown that for any quasihomomorphism $f\colon\mathbb{Z}\to\mathbb{Z}$
precisely one of the conditions
\begin{itemize}
\item $f$ has bounded range
\item for all $C>0$ there exists $n_{0}\in\mathbb{N}$ with $f(x)>C$ for
all $x>n_{0}$
\item for all $C>0$ there exists $n_{0}\in\mathbb{N}$ with $f(x)<-C$
for all $x>n_{0}$
\end{itemize}
holds. A real number $a=[f]$ is said to be positive if the second
condition holds for $f$. For all real numbers $b$ and $c$ it is
said that $b<c$ precisely when $c-b$ is positive.

\subsection{Knopfmacher-Knopfmacher's construction using Cantor's theorem (\cite{knopfmacher1987},
1987)\label{sub:KnopfmacherCantor}}

\subsubsection*{Observation (Cantor)}

Every real number $a>1$ can be written uniquely as 
\[
a=\prod_{k=1}^{\infty}(1+\frac{1}{a_{k}})=[a_{1},\dots,a_{k},\dots]
\]
where $a_{k}>0$ is an integer for all $k\ge1$ with $a_{k}>1$ from
some point onwards, and further $a_{k+1}\ge a_{k}^{2}$ for all $k\ge1$.
The number $a$ is rational if, and only if, $a_{k+1}=a_{k}^{2}$
for all $k\ge k_{0}$, for some $k_{0}\ge1$. Every real number $0<b<1$
can be written uniquely as 
\[
b=\prod_{k=1}^{\infty}(1-\frac{1}{b_{k}})=[-b_{1},\dots,-b_{k},\dots]
\]
where $b_{k}>1$ is an integer for all $k\ge1$ with $b_{k}>2$ from
some point onwards, and further $b_{k+1}>(b_{k}-1)^{2}$. The real
number $b$ is rational if, and only if, $b_{k+1}=1+(b_{k}-1)^{2}$
for all $k\ge k_{0}$, for some $k_{0}\ge1$.

\subsubsection*{The reals}

Let $S_{0}$ be the set of all sequences $[-b_{1},\dots,-b_{k},\dots]$
of negative integers $-b_{k}<-1$, not all equal to $-2$, and such
that $b_{k+1}>(b_{k}-1)^{2}$ for all $k\ge1$. Let $S_{1}$ be the
set of all sequences $[a_{1},\dots,a_{k},\dots]$ of positive integers
$a_{k}\ge1$, not all equal to $1$, and such that $a_{k+1}\ge a_{k}^{2}$
for all $k\ge1$. The set of non-negative real numbers is then $\mathbb{R}_{+}=S_{0}\cup\{1\}\cup S_{1}$.

\subsubsection*{Order}

For real numbers $a=[a_{1},\dots,a_{k},\dots],c=[c_{1},\dots,c_{k},\dots]\in S_{1}$
the relation $a<c$ holds precisely when for the first index $k_{0}$
where $a_{k_{0}}\ne c_{k_{0}}$ one has $a_{k_{0}}>b_{k_{0}}$. The
ordering among real numbers in $S_{0}$ is best seen by first introducing
the bijection $[a_{1},\dots,a_{k},\dots]\leftrightarrow[-b_{1},\dots,-b_{k},\dots]$
where $b_{k}=a_{k}+1$, which is in fact the reciprocal correspondence
$x\mapsto x^{-1}$. Then $a<c$ holds for real numbers $a,c\in S_{0}$
precisely when $c^{-1}>a^{-1}$ holds in $S_{1}$. Lastly, $a>1>b$
holds for all $a\in S_{1}$ and $b\in S_{0}$. The least upper bound
property for $\mathbb{R}$ is then proven.

\subsubsection*{Arithmetic}

The proof of Cantor's theorem yields an embedding of $\mathbb{Q}_{+}$
in $\mathbb{R}_{+}$ and further one obtains for every positive real
number $a$, by properly truncating it, sequences $a_{(n)}$ and $a^{(n)}$
in $\mathbb{R}_{+}$ of rational numbers which approximate $a$ from
below and from above, respectively. Addition of positive real numbers
$a$ and $b$ is then given by 
\[
a+b=\sup\{a_{(n)}+b_{(n)}\mid n\ge1\}
\]
and their product is given by 
\[
ab=\sup\{a_{(n)}b_{(n)}\mid n\ge1\}.
\]

\subsection{Pintilie's construction by infinite series (\cite{pintilie1988},
1988)\label{sub:Pintilie}}

\subsubsection*{Observation }

With the same starting point as in Shiu's construction described in
\subref{Shiu}, any sequence $(a_{n})$ of rational numbers with the
properties that $a_{n}\to0$ and $\sum a_{n}=\infty$ gives rise to
a presentation of the positive reals, i.e., for every real number
$a$ there exists a subsequence $(a_{n_{k}})$ with 
\[
a=\sum_{k=1}^{\infty}a_{n_{k}},
\]
albeit non-uniquely. In other words, if $A$ is the set of all subsequences
$(a_{n_{k}})$ such that 
\[
\sum_{k=1}^{\infty}a_{n_{k}}<\infty,
\]
then $\mathbb{R}_{+}\cong A$ as sets. To recover uniqueness one may
normalize the sequence $(a_{n})$ by demanding that $a_{0}=0$ and
only consider those subsequences leading to a bounded series which
further satisfy
\[
n_{k+1}=\min\{p>n_{k}\mid\exists m\colon\quad\sum_{i=1}^{k}a_{n_{i}}+a_{p}<\sum_{i=1}^{m}a_{n_{i}}\}
\]
for all $k\ge0$.

\subsubsection*{The reals}

Fix a sequence $(a_{n})_{n\ge1}$ of positive rational numbers and
set $a_{0}=0$. Let $A$ be the set of all subsequences $(a_{n_{k}})$
leading to a convergent series $\sum a_{n_{k}}$. Define the positive
real numbers to be the subset $\mathbb{R}_{+}\subseteq A$ consisting
only of those subsequences satisfying 
\[
n_{k+1}=\min\{p>n_{k}\mid\exists m\colon\quad\sum_{i=1}^{k}a_{n_{i}}+a_{p}<\sum_{i=1}^{m}a_{n_{i}}\}
\]
for all $k\ge0$.

\subsubsection*{Arithmetic}

Given positive real numbers $b=(b_{n})$ and $c=(c_{n})$, their sum
is given by the subsequence $(a_{n_{k}})$ determined, for all $k>0$,
by 
\[
n_{k+1}=\min\{n>n_{k}\mid\exists m\in\mathbb{N}\colon\quad\sum_{i=1}^{k}a_{n_{i}}+a_{n}<\sum_{i=1}^{m}(b_{i}+c_{i})\}
\]
and similarly their product is the subsequence determined by the conditions
\[
n_{k+1}=\min\{n>n_{k}\mid\exists m\in\mathbb{N}\colon\quad\sum_{i=1}^{k}a_{n_{i}}+a_{n}<(\sum_{i=1}^{m}b_{i})(\sum_{i=1}^{m}c_{i})\}
\]
for all $k>0$.

\subsubsection*{order}

For positive real numbers $b=(a_{b_{k}})$ and $c=(a_{c_{k}})$, the
meaning of $b<c$ is that there exists $k\in\mathbb{N}$ with $b_{k}>c_{k}$
and $b_{i}=c_{i}$, for all $i<k$.

\subsection{Knopfmacher-Knopfmacher's construction using Engel's theorem (\cite{knopfmacher1988},
1988)\label{sub:KnopfmacherEngel}}

\subsubsection*{Observation (Engel)}

Every real number $a$ can be written uniquely as 
\[
a=a_{0}+\frac{1}{a_{1}}+\frac{1}{a_{1}a_{2}}+\cdots+\frac{1}{a_{1}a_{2}\cdots a_{n}}+\cdots=(a_{0},a_{1},a_{2},\ldots)
\]
where the $a_{i}$ are integers satisfying $a_{i+1}\ge a_{i}\ge2$
for all $i\ge1$.

\subsubsection*{The reals}

Let $\mathbb{R}$ be the set of all infinite sequences $(a_{0},a_{1},a_{2},\dots)$
of integers satisfying $a_{k+1}\ge a_{k}\ge2$, for all $k\ge1$.

\subsubsection*{Order}

Given real numbers $A=(a_{0},a_{1},a_{2},\dots)$ and $B=(b_{0},b_{1},b_{2},\dots)$,
declare that $A<B$ precisely when 
\begin{itemize}
\item $a_{0}<b_{0}$, if $a_{0}\ne b_{0}$, or 
\item $a_{k}>b_{k}$ for the first index $k\ge1$ with $a_{k}\ne b_{k}$
otherwise.
\end{itemize}
The least upper bound property of $\mathbb{R}$ is then established
and the proof of Engel's theorem produces an order embedding $\mathbb{Q}\to\mathbb{R}$,
which thus serves to identify the rationals in $\mathbb{R}$. It then
follows that every real number can be approximated from above and
from below by, respectively, sequences $A_{(n)}$ and $A^{(n)}$ of
rationals.

\subsubsection*{Arithmetic}

Addition and multiplication of real numbers $A$ and $B$ is given
by exploiting the upper bound property of $\mathbb{R}$ and the rational
approximations above. That is, the sum $A+B$ is given by 
\[
A+B=\sup\{A_{(n)}+B_{(n)}\}
\]
and the product of positive reals is given by 
\[
AB=\sup\{A_{(n)}B_{(n)}\},
\]
and extended to all of $\mathbb{R}$ as usual. The proofs of the algebraic
properties utilize the rational approximations using limit-like arguments.

\subsection{Knopfmacher-Knopfmacher's construction using Sylvester's theorem
(\cite{knopfmacher1988}, 1988)}

The construction is formally identical to the one given in \subref{KnopfmacherEngel}
and will thus be presented quite briefly.

\subsubsection*{Observation (Sylvester)}

Every real number $a$ can be written uniquely as
\[
a=a_{0}+\frac{1}{a_{1}}+\frac{1}{a_{2}}+\cdots+\frac{1}{a_{n}}+\cdots=((a_{0},a_{1},a_{2},\dots))
\]
where the $a_{i}$ are integers satisfying $a_{1}\ge2$ and $a_{i+1}\ge a_{i}(a_{i}-1)+1$
for all $i\ge1$.

\subsubsection*{The reals}

Let $\mathbb{R}$ be the set of all infinite sequences $((a_{0},a_{1},a_{2},\dots))$
of integers satisfying $a_{k}\ge2$ and $a_{k+1}\ge a_{k}(a_{k}-1)+1$,
for all $k\ge1$.

\subsubsection*{Order}

Given real numbers $A=((a_{0},a_{1},a_{2},\dots))$ and $B=((b_{0},b_{1},b_{2},\dots))$,
declare that $A<B$ precisely when 
\begin{itemize}
\item $a_{0}<b_{0}$, if $a_{0}\ne b_{0}$, or 
\item $a_{k}>b_{k}$ for the first index $k\ge1$ with $a_{k}\ne b_{k}$
otherwise.
\end{itemize}
The least upper bound property of $\mathbb{R}$ is then established
and the proof of Sylvester's theorem produces an order embedding $\mathbb{Q}\to\mathbb{R}$,
identifying the rationals in $\mathbb{R}$. It then follows that every
real number can be approximated from above and from below by, respectively,
sequences $A_{(n)}$ and $A^{(n)}$ of rationals.

\subsubsection*{Arithmetic}

Formally identical to \subref{KnopfmacherEngel}
\begin{rem}
A generalization of Sylvester's theorem, and consequently a generalization
of this construction of the real numbers, is given, along essentially
the same lines, in \cite{SEL2014}. 
\end{rem}

\subsection{Knopfmacher-Knopfmacher's construction using the alternating Engel
theorem (\cite{knopfmacher1989}, 1989)\label{sub:AltEngel}}

\subsubsection*{Observation (alternating Engel)}

Every real number $A$ can be written uniquely as 
\[
A=a_{0}+\frac{1}{a_{1}}-\frac{1}{a_{1}a_{2}}+\frac{1}{a_{1}a_{2}a_{3}}-\cdots+\frac{(-1)^{n+1}}{a_{1}a_{2}\cdots a_{n}}+\cdots
\]
where the $a_{k}$ are integers satisfying $a_{k+1}\ge a_{k}+1\ge2$
for all $k\ge1$. Furthermore, this representation terminates after
a finite number of summands if, and only if, $A$ is rational.

\subsubsection*{The reals}

The set $\mathbb{R}$ of real numbers is defined to be the set of
all infinite sequences $(a_{0},a_{1},a_{2},\dots)$ of elements in
$\mathbb{N}^{+}$ which satisfy $a_{0}\in\mathbb{N}$ and $a_{k+1}\ge a_{k}+1\ge2$,
for all $k\ge1$.

\subsubsection*{Order}

Two real numbers $A=(a_{0},a_{1},a_{2},\dots)$ and $B=(b_{0},b_{1},b_{2},\dots)$
satisfy $A<B$ precisely when $a_{2n}<b_{2n}$ or $a_{2n+1}>b_{2n+1}$,
where the index $i=2n$ or $i=2n+1$ is the first index with $a_{i}\ne b_{i}$.
It is then shown that $\mathbb{R}$ satisfies the least upper bound
property. The proof of the alternating Engel theorem produces an order
embedding $\mathbb{Q}\to\mathbb{R}$, which thus serves to identify
the rationals in $\mathbb{R}$. It then follows that every real number
can be approximated from above and from below by, respectively, sequences
$A_{(n)}$ and $A^{(n)}$ of rationals.

\subsubsection*{Arithmetic}

Addition is defined by 
\[
A+B=\sup\{A^{(2n)}+B^{(2n)}\mid n\ge0\}
\]
and multiplication of positive reals is given by 
\[
AB=\sup\{A^{(2n)}B^{(2n)}\mid n\ge0\}.
\]
The rest of the construction is formally very similar to the one presented
in \subref{KnopfmacherEngel}.

\subsection{Knopfmacher-Knopfmacher's construction using the alternating Sylvester
theorem (\cite{knopfmacher1989}, 1989)}

\subsubsection*{Observation (alternating Sylvester)}

Every real number $A$ can be written uniquely as 
\[
A=a_{0}+\frac{1}{a_{1}}-\frac{1}{a_{2}}+\frac{1}{a_{3}}-\cdots+\frac{(-1)^{n=1}}{a_{n}}+\cdots=((a_{0},a_{1},a_{2},\ldots))
\]
where the $a_{k}$ are integers satisfying $a_{1}\ge1$ and $a_{k+1}\ge a_{k}(a_{k}+1)$
for all $k\ge1$. Furthermore, this representation terminates after
a finite number of summands if, and only if, $A$ is rational.

\subsubsection*{The reals}

The set $\mathbb{R}$ of real numbers is defined to be the set of
all infinite sequences $((a_{0},a_{1},a_{2},\dots))$ of elements
in $\mathbb{N}^{+}$ which satisfy, $a_{0}\in\mathbb{N}$, $a_{1}\ge1$,
and $a_{k+1}\ge a_{k}(a_{k}+1)$, for all $k\ge1$.

\subsubsection*{Order}

Formally identical to \subref{AltEngel}, except that it is the proof
of the alternating Sylvester theorem that gives the identification
of the rationals in $\mathbb{R}$.

\subsubsection*{Arithmetic}

Formally identical to \subref{AltEngel}.
\begin{rem}
In \cite{ikeda2013} a generalization of the alternating Sylvester
theorem is presented along with a corresponding construction of the
real numbers which is formally identical to this one. 
\end{rem}

\subsection{Arthan's irrational construction (\cite{Arthan}, 2001)}

\subsubsection*{Observation}

A closer look at the construction of $\mathbb{R}$ as a completion
of $\mathbb{Q}$ by means of Dedekind cuts (cf. \ref{sub:Dedekind})
reveals what the crucial ingredients present in $\mathbb{Q}$ actually
are that lead to the real numbers. In more detail, Dedekind's construction
is well-known to be a special case of the Dedekind-MacNeille completion
of an ordered set, and a famous theorem of Cantor shows that the completion
of any countable, unbounded, and densely ordered set is order isomorphic
to $\mathbb{R}$. Further, it is the archimedean property of $\mathbb{Q}$
that assures the additive structure on the completion has the right
properties, and finally, by a theorem of Hölder, the completion of
any ordered group which is dense and archimedean must be isomorphic
to an additive subgroup of $\mathbb{R}$, and thus admits a multiplication. 

It follows then that any countable, unbounded, archimedean, and densely
ordered group admits a completion isomorphic to $\mathbb{R}$ as a
field. If the multiplicative structure can effectively be defined
in terms of the given ordered group, then a construction of the reals
emerges.

\subsubsection*{The reals}

Given a dense, archimedean ordered commutative group, the reals are
constructed as its Dedekind-MacNeille completion. Once such an ordered
group is chosen the details are essentially identical to those of
Dedekind's construction and thus we further concentrate on the presentation
of the ordered group, namely $\mathbb{Z}[\sqrt{2}]$, which formally
we view as the set $\mathbb{Z}\times\mathbb{Z}$. Other rational numbers
may be chosen, with more or less adverse effects on the desired properties
of the group, the ease of establishing those properties, and the implementability
on a computer.

\subsubsection*{Arithmetic}

Addition and multiplication in $\mathbb{Z}[\sqrt{2}]$ are easily
established to be given by $(a,b)+(c,d)=(a+c,b+d)$ and $(a,b)(c,d)=(ac+2bd,ad+bc)$.

\subsubsection*{Order}

Recovering the ordering on $\mathbb{Z}[\sqrt{2}]$ solely in terms
of integers uses the fact that $\sqrt{2}$ is approximated by certain
rational numbers whose numerators and denominators admit an efficient
recursive formula. For details see \cite[Section 5.2]{Arthan}.

\subsection{Notes on Conway's surreal numbers and nonstandard models}

We conclude the survey by briefly mentioning two other venues leading
to the real numbers which, however, do not quite fall into the same
category as the constructions surveyed above.

\subsubsection{Conway's Surreal numbers (\cite{conway1976numbers}, 1976)}

Conway's famous construction of the surreal numbers is a construction
of a proper class in which every ordered field embeds. It thus follows
that the surreal numbers contain a copy of the real numbers and thus
one may view the surreal number system as providing yet another construction
of the real numbers. However, when one distills just the real numbers
from the entire array of surreal ones the construction basically collapses
to the Dedekind cuts construction. The interest in the surreal numbers
is not so much for the reals embedded in them, but rather for the
far reaching extra numbers beyond the reals that occupy most of the
surreal realm.

\subsubsection{Nonstandard constructions}

Since $\mathbb{R}$ is a completion of $\mathbb{Q}$ and since it
is well-known that the techniques of nonstandard analysis yield completions,
any path towards nonstandard analysis is also a path to a definition
of the real numbers. However, to what extent can these nonstandard
definitions be seen as \emph{constructions} of the real numbers is
a delicate issue. Inseparable to the technique of enlargement, which
is at the heart of nonstandard analysis, is the axiom of choice (or
some slightly weaker variant) and therefore the objects produced are
not particularly tangible. For that reason we avoided including any
details of nonstandard definitions of the real numbers. The interested
reader is referred to \cite{goldblatt1998lectures} for a very detailed
and carefully motivated exposition of one nonstandard model of the
reals (the hyperreals) and to \cite{benci2006eightfold} for seven
other possibilities. 

\bibliographystyle{plain}
\bibliography{generalReferences}

\end{document}